\documentclass[11 pt]{article}

\usepackage{amsmath,amsthm,amssymb,amsfonts,mathdots}
\usepackage{amsmath, amsfonts, amsthm, amssymb, mathtools, thmtools}
\usepackage[a4paper]{geometry}
\usepackage[ansinew]{inputenc}
\usepackage[numbers]{natbib}

\renewcommand{\Pr}{\mathrm{Pr}}

\DeclarePairedDelimiter{\abs}{\lvert}{\rvert}
\DeclarePairedDelimiter{\norm}{\lVert}{\rVert}

\newtheorem{theorem}{Theorem}

\newtheorem{corollary}{Corollary}

\newtheorem{definition}{Definition}

\newtheorem{proposition}{Proposition}

\begin{document}

\title{Adaptive Confidence Sets in $L^2$}

\author{Adam D. Bull ~~~~~~Richard Nickl \\   \\ University of Cambridge \footnote{Department of Pure Mathematics and Mathematical Statistics, Statistical Laboratory, Wilberforce Road CB30WB Cambridge, UK. Email: a.bull@statslab.cam.ac.uk, r.nickl@statslab.cam.ac.uk.}}

\date{\today ~~ [First version; October 7, 2011]}

\maketitle

\begin{abstract}
The problem of constructing confidence sets that are adaptive in $L^2$-loss over a continuous scale of Sobolev classes of probability densities is considered. Adaptation holds, where possible, with respect to both the radius of the Sobolev ball and its smoothness degree, and over maximal parameter spaces for which adaptation is possible. Two key regimes of parameter constellations are identified: one where full adaptation is possible, and one where adaptation requires critical regions be removed. Techniques used to derive these results include a general nonparametric minimax test for infinite-dimensional null- and alternative hypotheses, and new lower bounds for $L^2$-adaptive confidence sets.
\end{abstract}

%\tableofcontents

\section{Introduction}

The paradigm of adaptive nonparametric inference has developed a fairly complete theory for estimation and testing -- we mention the key references \cite{L90, DJKP95, DJKP96, LMS97, BBM99, BM01, S96} -- but the theory of adaptive confidence statements has not succeeded to the same extent, and consists in a significant part of negative results that are in a somewhat puzzling contrast to the fact that adaptive estimators exist. The topic of confidence sets is, however, of vital importance, since it addresses the question of whether the accuracy of adaptive estimation can itself be estimated, and to what extent the abundance of adaptive risk bounds and oracle inequalities in the literature are useful for statistical inference. 

In this article we give a set of necessary and sufficient conditions for when confidence sets that adapt to unknown smoothness in $L^2$-diameter exist in the problem of nonparametric density estimation on $[0,1]$. The scope of our techniques extends without difficulty to density estimation on the real line, and also to other common function estimation problems such as nonparametric regression or Gaussian white noise. Our focus on $L^2$-type confidence sets is motivated by the fact that they involve the most commonly used loss function in adaptive estimation problems, and so deserve special attention in the theory of adaptive inference. 

We can illustrate some main ideas by the simple example of two fixed Sobolev-type classes. Let $X_1, \dots, X_n$ be i.i.d.~with common probability density $f$ contained in the space $L^2$ of square-integrable functions on $[0,1]$. Let $\Sigma(r)=\Sigma(r,B)$ be a Sobolev ball of probability densities on $[0,1]$, of Sobolev-norm radius $B$ -- see Section \ref{inf} for precise definitions -- and consider adaptation to the submodel $\Sigma(s) \subset \Sigma(r)$, $s>r$. An adaptive estimator $\hat f_n$ exists, achieving the optimal rate $n^{-s/(2s+1)}$ for $f \in \Sigma(s)$ and $n^{-r/(2r+1)}$ otherwise, in $L^2$-risk; see for instance Theorem \ref{adapt} below. 

A confidence set is a random subset $C_n=C(X_1, \dots, X_n)$ of $L^2$. Define the $L^2$-diameter of a norm-bounded subset $C$ of $L^2$ as 
\begin{equation} \label{2diam}
|C| = \inf \left\{ \tau: C \subset \{h: \|h-g\|_2 \le \tau\} \text{ for some } g \in L^2 \right\},
\end{equation}
equal to the radius of the smallest $L^2$-ball containing $C$. For $G \subset L^2$ set $\|f-G\|_2= \inf_{g \in G}\|f-g\|_2$ and define, for $\rho_n \ge 0$ a sequence of real numbers, the separated sets $$\tilde \Sigma(r, \rho_n) \equiv \tilde \Sigma(r, s, B, \rho_n) = \{f \in \Sigma(r): \|f-\Sigma(s)\|_2 \ge \rho_n\}.$$ Obviously $\tilde \Sigma(r,0)=\Sigma(r)$, but for $\rho_n>0$ these sets are proper subsets of $\Sigma(r) \setminus \Sigma(s)$. We are interested in adaptive inference in the model $$\mathcal P_n \equiv \Sigma(s) \cup \tilde \Sigma(r, \rho_n)$$ under minimal assumptions on the size of $\rho_n$. We shall say that the confidence set $C_n$ is $L^2$-adaptive and honest for $\mathcal P_n$  if there exists a constant $M$ such that for every $n \in \mathbb N$,
\begin{equation}\label{fixad}
\sup_{f \in \Sigma(s)} {\Pr}_f\left\{|C_n| > M n^{-s/(2s+1)}\right\} \le \alpha',
\end{equation}
\begin{equation}\label{fixadr}
\sup_{f \in \tilde \Sigma(r,\rho_n)} {\Pr}_f\left\{|C_n| > M n^{-r/(2r+1)}\right\} \le \alpha'
\end{equation}
and if
\begin{equation} \label{hon}
\inf_{f \in \mathcal P_n} {\Pr}_f\left\{f \in C_n \right\} \ge 1-\alpha -r_n
\end{equation}
where $r_n \to 0$ as $n \to \infty$. We regard the constants $\alpha, \alpha'$ as given 'significance levels'. 

\begin{theorem} \label{fix}  Let $0<\alpha, \alpha ' <1, s>r>1/2$ and $B>1$ be given. \newline
A) An $L^2$-adaptive and honest confidence set for $\tilde \Sigma(r, \rho_n) \cup \Sigma(s)$ exists if one of the following conditions is satisfied: \newline
i) $s \le 2r$ and $\rho_n \ge 0$ \newline
ii) $s>2r$ and $$\rho_n \ge M n^{-r/(2r+1/2)}$$ for every $n \in \mathbb N$ and some constant $M$ that depends on $\alpha, \alpha', r, B$. \newline
B) If $s>2r$ and $C_n$ is an $L^2$-adaptive and honest confidence set for $\tilde \Sigma(r, \rho_n) \cup \Sigma(s)$, for every $\alpha, \alpha'>0$, then necessarily $$\liminf_n ~\rho_n n^{r/(2r+1/2)}>0.$$
\end{theorem}

We note first that for $s \le 2r$ adaptive confidence sets exist without any additional restrictions -- this is a main finding of the papers \cite{JL03, CL06, RV06} and has important precursors in \cite{L98, HL02, B04}. It is based on the idea that under the general assumption $f \in \Sigma(r)$ we may estimate the $L^2$-risk of any adaptive estimator of $f$ at precision $n^{-r/(2r+1/2)}$ which is $O(n^{-s/(2s+1)})$ precisely when $s \le 2r$. As soon as one wishes to adapt to smoothness $s>2r$, however, this cannot be used anymore, and adaptive confidence sets then require separation of $\Sigma(s)$ and $\Sigma(r) \setminus \Sigma(s)$ (i.e., $\rho_n>0$). Maximal subsets of $\Sigma(r)$ over which $L^2$-adaptive confidence sets do exist in the case $s>2r$ are given in Theorem \ref{fix}, with separation sequence $\rho_n$ characterised by the asymptotic order $n^{-r/(2r+1/2)}$.  This rate has, as we show in this article, a fundamental interpretation as the minimax rate of testing between the composite hypotheses 
\begin{equation} \label{compt}
H_0: f \in \Sigma(s) ~~ \text{against} ~~H_1: f \in \tilde \Sigma(r, \rho_n).
\end{equation} 
The occurrence of this rate in Theorem \ref{fix} parallels similar findings in Theorem 2 in \cite{HN11} in the different situation of confidence \textit{bands}, and is inspired by the general ideas in \cite{GN10, HN11, KNP11, B11}, which attempt to find 'maximal' subsets of the usual parameter spaces of adaptive estimation for which honest confidence statements can be constructed. Our results can be construed as saying that for $s>2r$ confidence sets that are $L^2$-adaptive exist precisely over those subsets of the parameter space $\Sigma(r)$ for which the target $s$ of adaptation is testable in a minimax way. 

Our solution of (\ref{compt}) is achieved in Proposition \ref{test} below, where we construct consistent tests for general composite problems of the kind $$H_0: f \in \Sigma ~~ \text{against} ~~H_1: f \in \Sigma(r), \|f-\Sigma\|_2 \ge \rho_n, ~~~\Sigma \subset \Sigma(r),$$ whenever the sequence $\rho_n$ is at least of the order $\max(n^{-r/(2r+1/2)}, r_n),$ where $r_n$ is related to the complexity of $\Sigma$ by an entropy condition. In the case $\Sigma = \Sigma(s)$ with $s>2r$ relevant here we can establish $r_n=n^{-s/(2s+1)}=o(n^{-r/(2r+1/2)}),$ so that this test is minimax in light of lower bounds in \cite{I86, I93}.

While the case of two fixed smoothness classes in Theorem \ref{fix} is appealing in its conceptual simplicity, it does not describe the typical adaptation problem, where one wants to adapt to a continuous smoothness parameter $s$ in a window $[r,R]$. Moreover the radius $B$ of $\Sigma(s)$ is, unlike in Theorem \ref{fix}, typically unknown, and the usual practise of 'undersmoothing' to deal with this problem incurs a rate-penalty for adaptation that we wish to avoid here.  Instead, we shall address the question of simultaneous exact adaptation to the radius $B$ and to the smoothness $s$. We first show that such strong adaptation is possible if $R<2r$, see Theorem \ref{rv}. In the general case $R\ge 2r$ we can use the ideas from Theorem \ref{fix} as follows: starting from a fixed largest model $\Sigma(r, B_0)$ with $r, B_0$ known, we discretise $[r,R]$ into a finite grid  $\mathcal S$ consisting of progressions $r, 2r, 4r, \dots$, and then use the minimax test for (\ref{compt}) in an iterated way to select the optimal value in $\mathcal S$. We then use the methods underlying Theorem \ref{fix} Ai) in the selected window, and show that this gives honest adaptive confidence sets over 'maximal' parameter subspaces $\mathcal P_n \subset\Sigma(r, B_0)$. In contrast to what is possible in the $L^\infty$-situation studied in \cite{B11}, the sets $\mathcal P_n$ asymptotically contain all of $\Sigma(r, B_0)$, highlighting yet another difference between the $L^2$- and $L^\infty$-theory.  See Proposition \ref{dish} and Theorem \ref{cont} below for details. We also present a new lower bound which implies that for $R>2r$ even 'pointwise in $f$' inference is impossible for the full parameter space of probability densities in the $r$-Sobolev space, see Theorem \ref{imp}. In other words, even asymptotically one has to remove certain subsets of the maximal parameter space if one wants to construct confidence sets that adapt to arbitrary smoothness degrees. One way to remove is to restrict the space apriori to a fixed ball $\Sigma(r,B_0)$ of known radius as discussed above, but other assumptions come to mind, such as 'self-similarity' conditions employed in \cite{PT00, GN10, KNP11, B11} for confidence intervals and bands. We discuss briefly how this applies in the $L^2$-setting.

We state all main results other than Theorem \ref{fix} above in Sections \ref{inf} and \ref{drei}, and proofs are given, in a unified way, in Section \ref{prf}

\section{The Setting} \label{inf}

\subsection{Wavelets and Sobolev-Besov Spaces} \label{sobsec}

Denote by $L^2:=L^2([0,1])$ the Lebesgue space of square integrable functions on $[0,1]$, normed by $\|\cdot\|_2$. For integer $s$ the classical Sobolev spaces are defined as the spaces of functions $f \in L^2$ whose (distributional) derivatives $D^\alpha f, 0<\alpha \le s,$ all lie in $L^2$. One can describe these spaces, for $s>0$ any real number, in terms of the natural sequence space isometry of $L^2$ under an orthonormal basis. We opt here to work with wavelet bases: for index sets $\mathcal Z \subset \mathbb Z, \mathcal Z_l \subset \mathbb Z$ and $J_0 \in \mathbb N$, let $$\{\phi_{J_0 m}, \psi_{lk}: m \in \mathcal Z, k \in \mathcal Z_l, l \ge J_0+1, l \in \mathbb N \}$$ be a compactly supported orthonormal wavelet basis of $L^2$ of regularity $S$, where as usual, $\psi_{lk}=2^{l/2}\psi_k(2^l\cdot)$. We shall only consider Cohen-Daubechies-Vial \cite{CDV93} wavelet bases where $|\mathcal Z_l|=2^l, |\mathcal Z| \le c(S)<\infty, J_0 \equiv J_0(S)$. We define, for $\langle f, g\rangle= \int_0^1 fg$ the usual $L^2$-inner product, and for $0 \le s <S$, the Sobolev (-type) norms
\begin{eqnarray} \label{sobolev}
\|f\|_{s,2}&:=& \max\left(2^{J_0s}\sqrt {\sum_{k \in \mathcal Z} \langle f, \phi_{J_0 k}\rangle^2},\sup_{l \ge J_0+1} 2^{ls} \sqrt{\sum_{k \in \mathcal Z_l}\langle f, \psi_{lk} \rangle^2 } \right) \notag \\
&=& \max \left(2^{J_0s}\|\langle f, \phi_{J_0 \cdot}\rangle\|_2,  \sup_{l \ge J_0+1} 2^{ls}\|\langle f, \psi_{l \cdot} \rangle\|_2 \right)
\end{eqnarray}
where in slight abuse of notation we use the symbol $\|\cdot\|_2$ for the sequence norms on $\ell^2(\mathcal Z_l), \ell^2(\mathcal Z)$ as well as for the usual norm on $L^2$. Define moreover the Sobolev (-type) spaces $$W^s \equiv B^s_{2\infty} = \{f \in L^2: \|f\|_{s,2}<\infty \}.$$ We note here that $W^s$ is not the classical Sobolev space -- in this case the supremum over $l \ge J_0+1$ would have to be replaced by summation over $l$ -- but the present definition gives rise to the slightly larger Besov space $B^s_{2 \infty}$, which will turn out to be the natural exhaustive class for our results below. We still refer to them as Sobolev spaces for simplicity, and since the main idea is to measure smoothness in $L^2$. We understand $W^s$ as spaces of continuous functions whenever $s>1/2$ (possible by standard embedding theorems). We shall moreover set, in abuse of notation, $\phi_{J_0 k} \equiv \psi_{J_0k}$ (which does not equal $2^{-1/2}\psi_{J_0+1,k}(2^{-1}\cdot)$) in order for the wavelet series of a function $f \in L^2$ to have the compact representation $$f=\sum_{l=J_0}^\infty \sum_{k \in \mathcal Z_{l}} \psi_{lk} \langle \psi_{lk},f\rangle,$$ with the understanding that $\mathcal Z_{J_0}=\mathcal Z$. The wavelet projection $\Pi_{V_j}(f)$ of $f \in L^2$ onto the span $V_j$ in $L^2$ of $$\{\phi_{J_0 m}, \psi_{lk}: m \in \mathcal Z, k \in \mathcal Z_l, J_0+1 \le l \le j\}$$ equals $$K_j(f)(x) \equiv \int_0^1 K_j(x,y)f(y)dy \equiv 2^j \int_0^1 K(2^jx, 2^jy)f(y)dy = \sum_{l=J_0}^{j-1} \sum_{k \in \mathcal Z_l}\langle f, \psi_{lk} \rangle \psi_{lk}(x) $$ where $K(x,y)=\sum_k \phi_{J_0 k}(x) \phi_{J_0 k}(y)$ is the wavelet projection kernel.

\subsection{Adaptive Estimation in $L^2$}

Let $X_1, \dots, X_n$ be i.i.d.~with common density $f$ on $[0,1]$, with joint distribution equal to the 
first $n$ coordinate projections of the infinite product probability measure ${\Pr}_f$. Write $E_f$ for the corresponding expectation operator. We shall throughout make the minimal assumption that $f \in W^r$ for some $r>1/2$, which implies in particular, by Sobolev's lemma, that $f$ is continuous and bounded on $[0,1]$. The adaptation problem arises from the hope that $f \in W^s$ for some $s$ significantly larger than $r$, without wanting to commit to a particular a priori value of $s$. In this generality the problem is still not meaningful, since the regularity of $f$ is not only described by containment in $W^s$, but also by the size of the Sobolev norm $\|f\|_{s,2}$. If one defines, for $0<s<\infty, 1 \le B<\infty$, the Sobolev-balls of densities 
\begin{equation} \label{ball}
\Sigma(s,B):= \left\{f:[0,1] \to [0, \infty), \int_T f =1, \|f\|_{s,2} \le B \right\},
\end{equation}
then Pinsker's minimax theorem (for density estimation) gives, as $n \to \infty$,
\begin{equation} \label{pinsker}
\inf_{T_n} \sup_{f \in \Sigma(s,B)} E_f \|T_n-f\|_2^2 \sim c(s) B^{2/(2s+1)} n^{-2s/(2s+1)}
\end{equation}
for some constant $c(s)>0$ depending only on $s$, and where the infimum extends over all measurable functions $T_n$ of $X_1, \dots, X_n$ (cf., e.g., the results in Theorem 5.1 in \cite{E08}). So any risk bound, attainable uniformly for elements $f \in \Sigma(s,B)$, cannot improve on $B^{2/(2s+1)}n^{-2s/(2s+1)}$ up to multiplicative constants. If $s,B$ are known then constructing estimators that attain this bound is possible, even with the asymptotically exact constant $c(s)$. The adaptation problem poses the question of whether estimators can attain such a risk bound without requiring knowledge of $B,s$. 

The paradigm of adaptive estimation has provided us with a positive answer to this problem, and one can prove the following result.  

\begin{theorem} \label{adapt}
Let $1/2 <r \le R<\infty$ be given. Then there exists an estimator $\hat f_n = f(X_1, \dots, X_n, r, R)$ such that, for every $s \in [r,R]$, every $B \ge 1, U>0$, and every $n \in \mathbb N$, $$\sup_{f \in \Sigma(s,B), \|f\|_\infty \le U} E_f \|\hat f_n - f\|_2^2 \le c B^{2/(2s+1)}n^{-2s/(2s+1)}$$
for a constant $0<c<\infty$ that depends only on $r, R, U$.
\end{theorem}

If one wishes to adapt to the radius $B \in [1,B_0]$ then the canonical choice for $U$ is 
\begin{equation} \label{supbd}
\sup_{f \in \Sigma(r,B_0)}\|f\|_\infty \le c(r) B_0 \equiv U < \infty,
\end{equation}
but other choices will be possible below. More elaborate techniques allow for $c$ to depend only on $s$, and even to obtain the exact asymptotic minimax 'Pinsker'-constant, see for instance Theorem 5.1 in \cite{E08}. We shall not study exact constants here, mostly to simplify the exposition and to focus on the main problem of confidence statements, but also since exact constants are asymptotic in nature and we prefer to give nonasymptotic bounds.

From a  'pointwise in $f$'  perspective we can conclude from Theorem \ref{adapt} that adaptive estimation is possible over the full continuous Sobolev scale $$\bigcup_{s \in [r,R], 1 \le B < \infty} \Sigma(s,B) = W^r \cap \left\{f:[0,1] \to [0, \infty), \int_0^1 f =1 \right\}; $$ for any probability density $f \in W^s, s \in [r,R]$, the single estimator $\hat f_n$ satisfies $$E_f\|\hat f_n -f\|_2^2 \le c \|f\|_{s,2}^{2/(2s+1)} n^{-2s/(2s+1)}$$ where $c$ depends on $r,R, \|f\|_\infty$. Since $\hat f_n$ does not depend on $B, U$ or $s$ we can say that $\hat f_n$ adapts to both $s \in [r,R]$ and $B \in [1, B_0]$ simultaneously. If one imposes an upper bound on $U$ then adaptation even holds for every $B \ge 1$. Our interest here is to understand what remains of this remarkable result if one is interested in adaptive \textit{confidence statements} rather than in risk bounds. 

\section{Adaptive Confidence Sets for Sobolev Classes} \label{drei}

\subsection{Honest Asymptotic Inference}

We aim to characterise those sets $\mathcal P_n$ consisting of uniformly bounded probability densities $f \in W^r$ for which we can construct adaptive confidence sets. More precisely, we seek random subsets $C_n$ of $L^2$ that depend only on known quantities, cover $f \in \mathcal P_n$ at least with prescribed probability $1-\alpha$, and have $L^2$-diameter $|C_n|$ adaptive with respect to radius and smoothness with prescribed probability at least $1-\alpha'$. To avoid discussing measurability issues we shall tacitly assume throughout that $C_n$ lies within an $L^2$-ball of radius $O(|C_n|)$ centered at a random variable $\tilde f_n \in L^2$. 

\begin{definition} [$L^2$-adaptive confidence sets]  \label{cc}
Let $X_1, \dots, X_n$ be i.i.d.~on $[0,1]$ with common density $f$. Let $0<\alpha, \alpha' <1$ and $1/2 <r \le R$ be given and let $C_n=C(X_1, \dots, X_n)$ be a random subset of $L^2$. $C_n$ is called $L^2$-adaptive and honest for a sequence of (nonempty) models $\mathcal P_n \subset W^r \cap \{f: \|f\|_\infty \le U\}$, if there exists a constant $L=L(r,R,U)$ such that for every $n \in \mathbb N$
\begin{equation} \label{adap}
\sup_{f \in \Sigma(s,B) \cap \mathcal P_n} {\Pr}_f\left\{|C_n| > L B^{1/(2s+1)}n^{-s/(2s+1)}\right\} \le \alpha'~~\text{for every}~s\in [r,R], B \ge 1,
\end{equation}
(the condition being void if $\Sigma(s,B) \cap \mathcal P_n$ is empty) and
\begin{equation} \label{hon}
\inf_{f \in \mathcal P_n} {\Pr}_f\left\{f \in C_n \right\} \ge 1-\alpha -r_n
\end{equation}
where $r_n \to 0$ as $n \to \infty$. 
\end{definition}

To understand the scope of this definition some discussion is necessary. First, the interval $[r,R]$ describes the range of smoothness parameters one wants to adapt to. Besides the restriction $1/2<r \le R < \infty$ the choice of this window of adaptation is arbitrary (although the values of $R,r$ influence the constants). Second, if we wish to adapt to $B$ in a fixed interval $[1,B_0]$ only, we may take $\mathcal P_n$ a subset of $\Sigma(r, B_0)$ and the canonical choice of $U=c(r)B_0$ from (\ref{supbd}). In such a situation (\ref{adap}) will still hold for every $B \ge 1$ although the result will not be meaningful for $B > B_0$. Otherwise we may impose an arbitrary uniform bound on $\|f\|_\infty$ and adapt to all $B \ge 1$. We require here the sharp dependence on $B$ in (\ref{adap}) and thus exclude the usual 'undersmoothed', near-adaptive, confidence sets in our setting.  A natural 'maximal' model choice would be $\mathcal P_n = \Sigma(r,B_0) ~ \forall n$ with $B_0 \ge 1$ arbitrary.

\subsection{The Case $R <  2r$.}
 
A first result, the key elements of which have been discovered and discussed in \cite{L98, HL02, JL03, CL06, RV06}, is that $L^2$-adaptive confidence statements that parallel the situation of Theorem \ref{adapt} exist without any additional restrictions whatsoever, in the case where $R < 2r$, so that the window of adaptation is $[r,2r)$. The sufficiency part of the following theorem is a simple extension of results in Robins and van der Vaart \cite{RV06} in that it shows that adaptation is possible not only to the smoothness $s$, but also to the radius $B$. The main idea of the proof is that, if $R<2r$, the squared $L^2$-risk of $\hat f_n$ from Theorem \ref{adapt} can be estimated at a rate compatible with adaptation, by a suitable $U$-statistic.

\begin{theorem} \label{rv} A) If $R<2r$, then for any $\alpha, \alpha'$, there exists a confidence set $C_n=C(X_1, \dots, X_n, r, R, \alpha, \alpha')$ which is honest and adaptive in the sense of Definition \ref{cc} for any choice $\mathcal P_n \equiv \Sigma(r,B_0) \cap \{f: \|f\|_\infty \le U\}, B_0 \ge 1, U>0$. \newline
B) If $R \ge 2r$, then for $\alpha, \alpha'$ small enough no $C_n$ as in A) exists. 
\end{theorem}

We emphasise that the confidence set $C_n$ constructed in the proof of Theorem \ref{rv} does only depend on $r,R,\alpha, \alpha'$ and does not require knowledge of $B_0$ or $U$. Note however that the sequence $r_n$ from Definition \ref{cc} does depend on $B_0$ -- one may thus use $C_n$ without any prior choice of parameters, but evaluation of its coverage is still relative to the model $\Sigma(r,B_0)$. Arbitrariness of $B_0, U$ implies, by taking  $B_0=\|f\|_{s,2}, U=\|f\|_\infty$ in the above result, that 'pointwise in $f$' adaptive inference is possible for any probability density in the Sobolev space $W^r$.

\begin{corollary} \label{rvp} 
Let $0<\alpha, \alpha' <1$ and $1/2 <r \le R$. Assume $R<2r$. There exists a confidence set $C_n=C(X_1, \dots, X_n, r, R, \alpha, \alpha')$ such that \newline
i) $\liminf_n {\Pr}_f\left\{f \in C_n \right\} \ge 1-\alpha~~~ \text{for every probability density } f \in  W^r,$ and
\newline
ii) $\limsup_n {\Pr}_f\{|C_n| > L \|f\|_{s, 2}^{1/(2s+1)}n^{-s/(2s+1)}\} \le \alpha'~~~ \text{for every probability density } f \in W^s, s\in [r,R],$ and some finite positive constant $L=L(r,R, \|f\|_\infty)$.
\end{corollary}

\subsection{The Case of General $R$ }

If we allow for general $R \ge 2r$ honest inference is not possible without restricting $\mathcal P_n$ further. In fact even a weaker 'pointwise in $f$' result of the kind of Corollary \ref{rvp} is impossible for general $R\ge r$. This is a consequence of the following lower bound.

\begin{theorem} \label{imp} 
Fix $0<\alpha<1/2$, let $s \ge r$ be arbitrary. A confidence set $C_n=C(X_1, \dots, X_n)$ in $L^2$ cannot satisfy \newline
i) $\liminf_n {\Pr}_f\{f \in C_n\} \ge 1- \alpha ~~~\text{for every probability density } f \in W^r$, and \newline
ii) $|C_n| = O_{{\Pr}_f}(r_n) ~~~\text{for every probability density } f \in W^s$ \newline
at any rate $r_n = o(n^{-r/(2r+1/2)})$.
\end{theorem} 

For $R>2r$ we have $n^{-R/(2R+1)} = o(n^{-r/(2r+1/2)})$. Thus even from a 'pointwise in $f$' perspective a confidence procedure cannot adapt to the entirety of densities in a Sobolev space $W^r$ when $R>2r$. On the other hand if we restrict to proper subsets of $W^r$, the situation may qualitatively change. For instance if we wish to adapt to submodels of a fixed Sobolev ball $\Sigma(r, B_0)$ with $r, B_0$ known, we have the following result.

\begin{proposition} \label{dish} Let $0<\alpha, \alpha' <1$ and $1/2 <r \le R, B_0 \ge 1$. There exists a confidence set  $C_n=C(X_1,\dots, X_n, B_0,r,R, \alpha, \alpha')$ such that \newline
i) $\liminf_n {\Pr}_f\left\{f \in C_n \right\} \ge 1-\alpha~~~ \text{for every probability density } f \in  \Sigma(r, B_0),$ and
\newline
ii) $\limsup_n {\Pr}_f\{|C_n| > L \|f\|_{s, 2}^{1/(2s+1)}n^{-s/(2s+1)}\} \le \alpha'~~~ \text{for every probability density } f \in \Sigma(s, B_0), s\in [r,R],$ and some finite positive constant $L=L(r,R, \|f\|_\infty)$.
\end{proposition}

Now if we compare Proposition \ref{dish} to Theorem \ref{rv} we see that there exists a genuine discrepancy between honest and pointwise in $f$ adaptive confidence sets when $R\ge 2r$. Of course Proposition \ref{dish} is not useful for statistical inference as the index $n$ from when onwards coverage holds depends on the unknown $f$. The question arises whether there are meaningful \textit{maximal} subsets of $\Sigma(r, B_0)$ for which honest inference is possible. The proof of Proposition \ref{dish} is in fact based on the construction of subsets $\mathcal P_n$ of $\Sigma(r,B_0)$ which grow dense in $\Sigma(r,B_0)$ and for which honest inference is possible. This approach follows the ideas from Part Aii) in Theorem \ref{fix}, and works as follows in the setting of continuous $s \in [r,R]$: assume without loss of generality that $2(N-1)r<R<2Nr$ for some $N \in \mathbb N, N >1$, and define the grid $$\mathcal S=\{s_m\}_{m=1}^N = \{r, 2r, 4r, \dots, 2(N-1)r\}.$$ Note that $\mathcal S$ is independent of $n$. Define, for $s \in \mathcal S \setminus \{s_N\}$,
\begin{equation*}
\tilde \Sigma(s, \rho):= \tilde \Sigma(s, B_0, \mathcal S, \rho) = \left\{f \in \Sigma(s, B_0): \|f-\Sigma(t, B_0)\|_2 \ge \rho ~\forall t>s, t \in \mathcal S \right\}.
\end{equation*}
We will choose the separation rates $$\rho_n(s) \sim n^{-s/(2s+1/2)},$$ equal to the minimax rate of testing between $\Sigma(s, B_0)$ and any submodel $\Sigma(t, B_0)$ for $t \in \mathcal S, t>s$. The resulting model is therefore, for $M$ some positive constant, $$\mathcal P_n(M, \mathcal S) = \Sigma(s_N, B_0) \bigcup \left(\bigcup_{s \in \mathcal S \setminus \{s_N\}} \tilde \Sigma(s, M\rho_n(s))\right).$$ 

The main idea behind the following theorem is to first construct a minimax test for the nested hypotheses $$\{H_s: f \in \tilde \Sigma(s, M \rho_n(s))\}_{s \in \mathcal S \setminus \{s_N\}},$$ then to estimate the risk of the adaptive estimator $\hat f_n$ from Theorem \ref{adapt} under the assumption that $f$ belongs to smoothness hypothesis selected by the test, and to finally construct a confidence set centered at $\hat f_n$ based on this risk estimate (as in the proof of Theorem \ref{rv}).

\begin{theorem} \label{cont}  Let $R > 2r$ and $B_0 \ge 1$ be arbitrary. There exists a confidence set $C_n=C(X_1,\dots, X_n, B_0,r,R, \alpha, \alpha')$, honest and adaptive in the sense of Definition \ref{cc}, for $\mathcal P_n = \mathcal P_n(M, \mathcal S), n \in \mathbb N,$ with $M$ a large enough constant and $U$ as in (\ref{supbd}).
\end{theorem}

First note that, since $\mathcal S$ is independent of $n$, $\mathcal P_n(M, \mathcal S) \nearrow \Sigma(r, B_0)$ as $n \to \infty$, so that the model $\mathcal P_n(M, \mathcal S)$ grows dense in the fixed Sobolev ball, which for known $B_0$ is the full model. This implies in particular Proposition \ref{dish}. 

An important question is whether $\mathcal P_n(M, \mathcal S)$ was taken to grow as fast as possible as a function of $n$, or in other words, whether a smaller choice of $\rho_n(s)$ would have been possible.  The lower bound in Theorem \ref{fix} implies that any faster choice for $\rho_n(s)$ makes honest inference impossible. Indeed, if $C_n$ is an honest confidence set over $\mathcal P_n(M, \mathcal S)$ with a faster separation rate $\rho_n'=o(\rho_n(s))$ for some $s \in \mathcal S \setminus \{s_N\}$, then we can use $C_n$ to test $H_0: f \in \Sigma(s')$ against $H_1: f \in \tilde \Sigma(s, \rho_n')$ for some $s'>2s$, which by the proof of Theorem \ref{fix} gives a contradiction.

\subsubsection{Self-Similarity Conditions} \label{ssc}

The proof of Theorem \ref{cont} via testing smoothness hypotheses is strongly tied to knowledge of the upper bound $B_0$ for the radius of the Sobolev ball, but as discussed above, this cannot be avoided without contradicting Theorem \ref{imp}. Alternative ways to restrict \(W^r,\) other than constraining the radius, and which may be practically relevant, 
are given in \cite{PT00, GN10, KNP11, B11}. The authors instead restrict to `self-similar' functions, whose 
regularity is similar at large and small scales. As the results \cite{GN10, KNP11, B11} prove adaptation in \(L^\infty,\) they naturally imply adaptation also in \(L^2;\) the functions excluded, however, are now those whose norm is hard to estimate, rather than those whose norm is merely large. In the $L^2$-case we need to estimate $s$ only up to a small constant; as this is more favourable than the $L^\infty$-situation, one may impose weaker self-similarity assumptions, tailored to the $L^2$-situation. This can be achieved arguing in a similar fashion to Bull \cite{B11}, but we do not pursue this further in the present paper.

\section{Proofs} \label{prf}

\subsection{Some Concentration Inequalities} \label{ci}

Let $X_i, i=1, 2, \dots,$ be the coordinates of the product probability space $(T,{\cal T},P)^{\mathbb N}$, where $P$ is any probability measure on $(T, \mathcal T)$, $P_n=n^{-1}\sum_{i=1}^n \delta_{X_i}$ the empirical measure, $E$ expectation under $P^\mathbb N \equiv \Pr$. For $M$ any set and $H:M \to \mathbb R$, set $\|H\|_M=\sup_{m \in M}|H(m)|$. We also write $Pf=\int_T fdP$ for measurable $f: T \to \mathbb R$.

The following Bernstein-type inequality for canonical $U$-statistics of order two is due to Gin\'{e}, Latala and Zinn \cite{GLZ00}, with refinements about the numerical constants in Houdr\'{e} and Reynaud-Bouret \cite{HR03}: let $R(x,y)$ be a symmetric real-valued function defined on $T \times T$, such that $ER(X,x)=0$ for all $x$, and let
$$\Lambda^2_1=\frac{n(n-1)}{2} ER(X_1,X_2)^2,$$
$$\Lambda_2=n\sup\{E[R(X_1,X_2)\zeta(X_1)\xi(X_2)]:E\zeta^2(X_1)\le 1,E\xi^2(X_1)\le1\},$$
$$ \Lambda_3=\|nER^2(X_1,\cdot)\|^{1/2}_\infty,\ \ \Lambda_4=\|R\|_\infty.$$
Let moreover $U_n^{(2)}(R) = \frac{2}{n(n-1)} \sum_{i<j} R(X_i, X_j)$ be the corresponding degenerate $U$-statistic of order two. Then, there exists a universal constant $0<C<\infty$ such that for all $u>0$ and $n \in\mathbb N$:
\begin{equation} \label{glz}
\Pr\left\{\frac{n(n-1)}{2}|U_n^{(2)}(R)|>C(\Lambda_1u^{1/2}+\Lambda_2u+\Lambda_3u^{3/2}+\Lambda_4u^2)\right \} \le 6\exp\{- u\}.
\end{equation}

We will also need Talagrand's \cite{T96} inequality for empirical processes. Let $\cal F$ be a countable class of measurable functions on $T$ that take values in $[-1/2,1/2]$, or, if $\mathcal F$ is $P$-centered, in $[-1,1]$. Let $\sigma \le 1/2$, or $\sigma \le 1$ if $\mathcal F$ is $P$-centered, and $V$ be any two numbers satisfying
\begin{equation*}
\sigma^2 \geq \|Pf^2\|_{\cal F},\ \ V \geq n\sigma^2+2E\left\|\sum_{i=1}^n(f(X_i)-Pf)\right\|_{\cal F}.
\end{equation*}
Bousquet's \cite{B03} version of Talagrand's inequality then states: for every $u >0$,
\begin{equation}\label{bous}
\Pr\left\{\left\|\sum_{i=1}^n(f(X_i)-Pf)\right\|_{\cal F}\ge E\left\|\sum_{i=1}^n(f(X_i)-Pf)\right\|_{\cal F}+u\right\}\le \exp\left(-\frac{u^2}{2V+\frac{2}{3}u}\right).
\end{equation}
A consequence of this inequality, derived in Section 3.1 in \cite{GN11}, is the following. If $T=[0,1]$, $P$ has bounded Lebesgue density $f$ on $T$, and $f_n(j)=\int_0^1 K_j(\cdot,y)dP_n(y)$, then for $M$ large enough, every $j \ge 0, n \in \mathbb N$ and some positive constants $c, c'$ depending on $U$ and the wavelet regularity $S$,
\begin{equation} \label{tal2}
\sup_{f: \|f\|_\infty \le U}{\Pr}_f \left \{ \left\|f_n(j) - Ef_n(j) \right\|_2 > M  \sqrt{\|f\|_\infty \frac{2^j}{n}} \right \} \le c' e^{-cM^2 2^j}.
\end{equation}

\subsection{A General Purpose Test for Composite Nonparametric Hypotheses} \label{testsec}

In this subsection we construct a general test for composite nonparametric null hypotheses that lie in a fixed Sobolev ball, under assumptions only on the entropy of the null-model. While of independent interest, the result will be a key step in the proofs of Theorems \ref{fix} and \ref{cont}.

Let $X,X_1, \dots, X_n$ be i.i.d.~with common probability density $f$ on $[0,1]$, let $\Sigma$ be any subset of a fixed Sobolev ball $\Sigma (t,B)$ for some $t>1/2$ and consider testing
\begin{equation} \label{genhyp}
H_0: f \in \Sigma~\text{ against } H_1: f \in \Sigma(t, B)\setminus \Sigma, \|f-\Sigma\|_2 \ge \rho_n ,
\end{equation}
where $\rho_n \ge 0$ is a sequence of nonnegative real numbers. For $\{\psi_{lk}\}$ a $S$-regular wavelet basis, $S>t$, $J_n \ge J_0$ a sequence of positive integers such that $2^{J_n} \simeq n^{1/(2t+1/2)}$ and for $g \in \Sigma$, define the $U$-statistic
\begin{equation}
T_n(g) = \frac{2}{n (n-1)} \sum_{i<j} \sum_{l=J_0}^{J_n-1} \sum_{k \in \mathcal Z_l} (\psi_{lk}(X_i)-\langle \psi_{lk}, g\rangle)(\psi_{lk}(X_j)-\langle \psi_{lk}, g \rangle)
\end{equation}
and, for $\tau_n$ some thresholds to be chosen below, the test statistic
\begin{equation} \label{stat}
\Psi_n = 1\left\{\inf_{g \in \Sigma} |T_n(g)| > \tau_n \right\}.
\end{equation}
Measurability of the infimum in (\ref{stat}) can be established by standard compactness/continuity arguments.

We shall prove a bound on the sum of the type-one and type-two errors of this test under some entropy conditions on $\Sigma$, more precisely, on the class of functions $$\mathcal G(\Sigma) = \bigcup_{J > J_0} \left\{\sum_{l=J_0}^{J-1} \sum_{k \in \mathcal Z_l} \psi_{lk}(\cdot) \langle \psi_{lk}, g \rangle: g \in \Sigma \right\}.$$  Recall the usual covering numbers $N(\varepsilon, \mathcal G, L^2(P))$ and bracketing metric entropy numbers $N_{[]}(\varepsilon, \mathcal G, L^2(P))$ for classes $\mathcal G$ of functions and probability measures $P$ on $[0,1]$ (e.g., \cite{G00, VW96}).

\begin{definition} \label{entropy}
Say that $\Sigma$ is $s$-regular if one of the following conditions is satisfied for some fixed finite constants $A$ and every $0<\varepsilon <A$:
\newline a) For any probability measure $Q$ on $[0,1]$ (and $A$ independent of $Q$) we have $$\log N(\varepsilon, \mathcal G(\Sigma), L^2(Q)) \le (A/\varepsilon)^{1/s}.$$
\newline
b) For $P$ such that $dP=fd\lambda$ with Lebesgue density $f:[0,1] \to [0, \infty)$ we have $$\log N_{[]}(\varepsilon, \mathcal G(\Sigma), L^2(P)) \le (A/\varepsilon)^{1/s}.$$ 
\end{definition}

Note that a ball $\Sigma(s,B)$ satisfies this condition for the given $s, 1/2<s<S,$ since any element of $\mathcal G(\Sigma(s,B))$ has $\|\cdot\|_{s,2}$-norm no more than $B$, and since $$\log N(\varepsilon, \Sigma(s,B), \|\cdot\|_\infty) \le (A/\varepsilon)^{1/s},$$ see, e.g., p.506 in \cite{LGM96}.

\begin{proposition} \label{test} Let $$\tau_n = L d_n \max(n^{-2s/(2s+1)}, n^{-2t/(2t+1/2)}), ~~~\rho^2_n = \frac{L_0}{L} \tau_n$$ for real numbers $1 \le d_n \le d(\log n)^\gamma$ and positive constants $L, L_0, \gamma,d$. Let the hypotheses $H_0, H_1$ be as in (\ref{genhyp}), the test $\Psi_n$ as in (\ref{stat}), and assume $\Sigma$ is $s$-regular for some $s>1/2$.  Then for $L=L(B, t, S)$, $L_0=L_0(L, B,t, S)$ large enough and every $n \in \mathbb N$ there exist constants $c_i, i=1,\dots, 3$ depending only on $L,L_0, t, B$ such that $$\sup_{f \in H_0} E_f\Psi_n + \sup_{f \in H_1} E_f(1-\Psi_n) \le c_1 e^{-d_n^2} +  c_2 e^{-c_3 n \rho_n^2}.$$ 
\end{proposition} 
The main idea of the proof is as follows: for the type-one errors our test-statistic is dominated by a degenerate $U$-statistic which we can bound with inequality (\ref{glz}), carefully controlling the four regimes present. For the alternatives the test statistic can be decomposed into a degenerate $U$-statistic which can be dealt with as before, and a linear part, which is the critical one. The latter can be compared to a ratio-type empirical process which we control by a slicing argument applied to $\Sigma$, combined with Talagrand's inequality.

\begin{proof} 1) We first control the type-one errors. Since $f \in H_0 = \Sigma$ we see
\begin{equation} \label{h0}
E_f\Psi_n = {\Pr}_f \left\{\inf_{g \in \Sigma} |T_n(g)| > \tau_n \right\} \le {\Pr}_f \left\{|T_n(f)| > \tau_n \right\}.
\end{equation}
$T_n(f)$ is a $U$-statistic with kernel $$R_f(x,y)= \sum_{l=J_0}^{J_n-1} \sum_{k \in \mathcal Z_l}(\psi_{lk}(x)-\langle \psi_{lk}, f\rangle)(\psi_{lk}(y)-\langle \psi_{lk}, f \rangle),$$ which satisfies $E R_f(x,X_1)=0$ for every $x$, since $E_f(\psi_{lk}(X)-\langle \psi_{lk}, f\rangle)=0$ for every $k,l$. Consequently $T_n(f)$ is a degenerate $U$-statistic of order two, and we can apply inequality (\ref{glz}) to it, which we shall do with $u=d^2_n$. We thus need to bound the constants $\Lambda_1, \dots, \Lambda_4$ occurring in inequality (\ref{glz}) in such a way that, for $L$ large enough,

\begin{equation} \label{bd0}
\frac{2C}{n(n-1)}(\Lambda_1 d_n+\Lambda_2 d_n^2+\Lambda_3 d_n^3+\Lambda_4 d_n^4) \le L d_n n^{-2t/(2t+1/2)} \le \tau_n,
\end{equation}
which is achieved by the following estimates, noting that $n^{-2t/(2t+1/2)} \simeq 2^{J_n/2}/n$. 

First, by standard $U$-statistic arguments, we can bound $ER^2_f(X_1,X_2)$ by the second moment of the uncentred kernel, and thus, using orthonormality of $\psi_{lk}$,
\begin{eqnarray*}
ER_f^2(X_1,X_2) &\le& \int \int \left(\sum_{k,l} \psi_{lk}(x) \psi_{lk}(y)\right)^2 f(x)f(y)dxdy  \\
&\le& \|f\|_\infty ^2 \sum_{l=J_0}^{J_n-1} \sum_{k \in \mathcal Z_l} \int_0^1 \psi_{lk}^2(x)dx  \int_0^1 \psi_{lk}^2(y)dy  \\
&\le & C(S)2^{J_n} \|f\|^2_\infty
\end{eqnarray*}
for some constant $C(S)$ that depends only on the wavelet basis. We obtain $\Lambda_1^2 \leq C(S) n(n-1) 2^{J_n} \|f\|_\infty^2/2$ and it follows, using (\ref{supbd}) that for $L$ large enough and every $n$, $$\frac{2C\Lambda_1 d_n}{n(n-1)} \le C(S, B, t) \frac{2^{J_n/2}d_n}{n} \le \tau_n/4.$$ For the second term note that, using the Cauchy-Schwarz inequality and that $K_j$ is a projection operator
\begin{eqnarray*} 
\left|\int \int \sum_{l=J_0}^{J_n-1} \sum_{k \in \mathcal Z_l}\psi_{lk}(x) \psi_{lk}(y) \zeta(x) \xi(y) f(x)f(y)dxdy \right| &=& \left|\int K_{J_n} (\zeta f)(y) \xi(y) f(y)dy \right| \\
&\le& \|K_{J_n}(\zeta f)\|_2 \|\xi f\|_2 \le \|f\|_\infty^2,
\end{eqnarray*}
and similarly 
\begin{equation*}
|E[E_{X_1} [K_{J_n}(X_1,X_2)] \zeta(X_1) \xi(X_2)]| \leq \|f\|^2_\infty, \ |EK_{J_n}(X_1,X_2)| \leq \|f\|^2_\infty.
\end{equation*}
Thus
$$E[R_f(X_1,X_2)\zeta(X_1)\xi(X_2)] \leq 4\|f\|^2_\infty $$
so that, using (\ref{supbd}), $$\frac{2C\Lambda_2 d_n^2}{n(n-1)} \le \frac{C'(B,t) d_n^2}{n} \le \tau_n/4$$ again for $L$ large enough and every $n$.  

For the third term, using the decomposition $R_f(x_1,x)=(r(x_1,x)-E_{X_1}r(X,x))+(E_{X,Y}r(X,Y)-E_Yr(x_1,Y))$ for $r(x,y)=\sum_{k,l} \psi_{lk}(x)\psi_{lk}(y)$, the inequality $(a+b)^2 \le 2a^2+2b^2$ and again orthonormality, we have that for every $x\in\mathbb R$, 
$$n|E_{X_1}R_f^2(X_1,x)| \leq 2n \left[\|f\|_\infty \sum_{l=J_0}^{J_n-1}\sum_{k \in \mathcal Z_l} \psi^2_{lk}(x) + \|f\|_\infty \|\Pi_{V_{J_n}}(f)\|_2^2 \right]$$ 
so that, using $\|\psi_{lk}\|_\infty \le d2^{l/2}$, again for $L$ large enough and by (\ref{supbd}),
$$\frac{2C\Lambda_3d_n^{3}}{n (n-1)}\le C''(B,t) \frac{2^{J_n/2}d_n^3}{n}\frac{1}{\sqrt n} \le \tau_n/4.$$
Finally, we have $\Lambda_4=\|R_f\|_\infty \leq c 2^{J_n}$ and hence $$\frac{2C\Lambda_4 d_n^4}{n(n-1)} \le C' \frac{2^{J_n} d_n^4}{n^2} \le \tau_n/4,$$ so that we conclude for $L$ large enough and every $n \in \mathbb N$, from inequality (\ref{glz}),
\begin{equation}
{\Pr}_f \left\{|T_n(f)| > \tau_n \right\} \le 6\exp\left\{-d^2_n\right\}
\label{second}
\end{equation}
which completes the bound for the type-one errors in view of (\ref{h0}).

2) We now turn to the type-two errors. In this case, for $f \in H_1$
\begin{equation} \label{h1}
E_f(1-\Psi_n) = {\Pr}_f \left\{ \inf_{g \in \Sigma} |T_n(g)| \le \tau_n \right\}.\end{equation}
and the typical summand of $T_n(g)$ has Hoeffding-decomposition
\begin{eqnarray*}
&& (\psi_{lk}(X_i)-\langle \psi_{lk},g \rangle)(\psi_{lk}(X_j)-\langle \psi_{lk}, g \rangle) \\
&& = (\psi_{lk}(X_i)-\langle \psi_{lk}, f \rangle + \langle \psi_{lk},f-g\rangle)(\psi_{lk}(X_j)-\langle \psi_{lk}, f \rangle + \langle \psi_{lk}, f-g \rangle) \\
&& = (\psi_{lk}(X_i)- \langle \psi_{lk}, f \rangle)(\psi_{lk}(X_j)- \langle \psi_{lk}, f \rangle)) \\
&& ~~~~+ (\psi_{lk}(X_i)-\langle \psi_{lk}, f\rangle) \langle \psi_{lk}, f-g \rangle + (\psi_{lk}(X_j)-\langle \psi_{lk}, f \rangle) \langle \psi_{lk}, f-g \rangle \\
&& ~~~~+ \langle \psi_{lk}, f-g \rangle^2
\end{eqnarray*}
so that by the triangle inequality, writing 
\begin{equation} \label{linear}
L_n(g)= \frac{2}{n} \sum_{i=1}^n  \sum_{l=J_0}^{J_n-1}\sum_{k \in \mathcal Z_l}(\psi_{lk}(X_i)-\langle \psi_{lk}, f\rangle) \langle \psi_{lk}, f-g \rangle
\end{equation}
for the linear terms, we conclude
\begin{eqnarray} \label{2lb}
\left| T_n(g) \right | &\ge& \sum_{l=J_0}^{J_n-1}\sum_{k \in \mathcal Z_l}\langle \psi_{lk}, f-g \rangle^2 - \left|T_n(f)\right| -|L_n(g)| \notag \\
& = & \|\Pi_{V_{J_n}}(f-g)\|_2^2 - |T_n(f)| - |L_n(g)|
\end{eqnarray}
for every $g \in \Sigma$. 

We can find random $g^*_n \in \Sigma$ such that $\inf_{g \in \Sigma} |T_n(g)| = |T_n(g^*_n)|$. (If the infimum is not attained the proof below requires obvious modifications; for the case $\Sigma=\Sigma(s,B), s>t$, relevant below, the infimum can be shown to be attained at a measurable minimiser by standard continuity and compactness arguments.) We bound the probability in (\ref{h1}), using (\ref{2lb}), by $$ {\Pr}_f \left\{|L_n(g_n^*)| >  \frac{\|\Pi_{V_{J_n}}(f-g_n^*)\|_2^2 - \tau_n}{2}\right\} +{\Pr}_f \left\{|T_n(f)| > \frac{\|\Pi_{V_{J_n}}(f-g_n^*)\|_2^2-\tau_n}{2}\right\}.$$ Now by the standard approximation bound (cf.~(\ref{sobolev})) and since $g^*_n \in \Sigma \subset \Sigma(t, B)$, 
\begin{equation} \label{2sep}
\|\Pi_{V_{J_n}}(f-g_n^*)\|_2^2 \ge \inf_{g \in \Sigma}\|f-g\|_2^2 - c(B)2^{-2J_nt} \ge 4 \tau_n
\end{equation}
 for $L_0$ large enough depending only on $B$ and the choice of $L$ from above. We can thus bound the sum of the last two probabilities by
 $$ {\Pr}_f \{|L_n(g_n^*)| >  \|\Pi_{V_{J_n}}(f-g_n^*)\|_2^2/4\} +{\Pr}_f \{|T_n(f)| > \tau_n\}.$$ For the second degenerate part the proof of Step 1 applies, as only boundedness of $f$ was used there. In the linear part somewhat more care is necessary. We have
\begin{eqnarray} \label{ratio}
{\Pr}_f \{|L_n(g^*_n)| > \|\Pi_{V_{J_n}}(f-g^*_n)\|_2^2/4\} \le {\Pr}_f \left\{\sup_{g \in \Sigma}\frac{|L_n(g)|}{\|\Pi_{V_{J_n}}(f-g)\|_2^2} > \frac{1}{4}\right\}.
\end{eqnarray}
Note that the variance of the linear process from (\ref{linear}) can be bounded, for fixed $g \in \Sigma$, using independence and orthonormality, by
\begin{eqnarray} \label{weakvar}
Var_f(|L_n(g)|) &\le& \frac{4}{n} \int \left(\sum_{l=J_0}^{J_n-1}\sum_{k \in \mathcal Z_l} \psi_{lk}(x) \langle \psi_{lk}, f-g\rangle\right)^2 f(x)dx \notag \\
&\le& \frac{4 \|f\|_\infty}{n} \sum_{l=J_0}^{J_n-1}\sum_{k \in \mathcal Z_l} \int \psi^2_{lk}(x)dx \cdot \langle \psi_{lk}, f-g \rangle^2 \notag \\
&\le& \frac{4\|f\|_\infty \|\Pi_{V_{J_n}}(f-g)\|_2^2}{n}
\end{eqnarray}
so that the supremum in (\ref{ratio}) is one of a self-normalised ratio-type empirical process. Such processes can  be controlled by slicing the supremum into shells of almost constant variance, cf.~Section 5 in \cite{G00} or \cite{GK06}. Define, for $g \in \Sigma$, $$\sigma^2(g):=\|\pi_{V_{J_n}}(f-g)\|_2^2 \ge \|f-g\|_2^2 - c(B)2^{-2J_n t} \ge c \rho_n^2,$$ the inequality holding for $L_0$ large enough and some $c>0$, as in (\ref{2sep}). Define moreover, for $m \in \mathbb Z$, the class of functions $$\mathcal G_{m, J_n}  = \left\{2\sum_{l=J_0}^{J_n-1}\sum_{k \in \mathcal Z_l} \psi_{lk}(\cdot) \langle \psi_{lk}, f-g \rangle: g \in \Sigma, \sigma^2(g)\le 2^{m+1} \right\},$$ which is uniformly bounded by a constant multiple of $\|f\|_{t,2}+\sup_{g \in \Sigma(t,B)}\|g\|_{t,2} \le 2B$ in view of (\ref{sobolev}) and since $t>1/2$. Then clearly, in the notation of Subsection \ref{ci}, $$\sup_{g \in \Sigma: \sigma^2(g) \le 2^{m+1}}|L_n(g)| = \|P_n-P\|_{\mathcal G_{m, J_n}} $$ and we bound the last probability in (\ref{ratio}) by
\begin{eqnarray} \label{sup}
&& {\Pr}_f \left\{\max_{m \in \mathbb Z: c'\rho_n^2 \le 2^m \le C}\sup_{g \in \Sigma: 2^m \le \sigma^2(g) \le 2^{m+1}}\frac{|L_n(g)|}{\sigma^2(g)} > \frac{1}{4}\right\} \notag \\
&& \le \sum_{m \in \mathbb Z: c'\rho_n^2 \le 2^m \le C} {\Pr}_f \left\{\sup_{g \in \Sigma: \sigma^2(g) \le 2^{m+1}}|L_n(g)| > 2^{m-2}\right\} \\
&& \le  \sum_{m \in \mathbb Z: c'\rho_n^2 \le 2^m \le C} {\Pr}_f \left\{ \|P_n-P\|_{\mathcal G_{m, J_n}}-E \|P_n-P\|_{\mathcal G_{m, J_n}} > 2^{m-2} -E \|P_n-P\|_{\mathcal G_{m, J_n}}\right\} \notag
\end{eqnarray}
where we may take $C<\infty$ as $\Sigma \subset \Sigma(t, B)$ is bounded in $L^2$, and where $c'$ is a positive constant such that $c' \rho_n^2 \le 2^m \le c \rho_n^2$ for some $m \in \mathbb Z$. We bound the expectation of the empirical process. Both the uniform and the bracketing entropy condition for $\mathcal G(\Sigma)$ carry over to $\cup_{J \ge 0}\mathcal G_{J,m}$ since translation by $f$ preserves the entropy. Using the standard entropy-bound plus chaining moment inequality (3.5) in Theorem 3.1 in \cite{GK06} in case a) of Definition \ref{entropy}, and the second bracketing entropy moment inequality in Theorem 2.14.2 in \cite{VW96} in case b), together with the variance bound (\ref{weakvar}) and with (\ref{supbd}), we deduce 
\begin{equation} \label{mom}
E \|P_n-P\|_{\mathcal G_{m,J_n}} \le  C \left( \sqrt{\frac{2^m}{n}} (2^m)^{-1/4s} + \frac{(2^m)^{-1/2s}}{n}\right).
\end{equation}
We see that $$2^{m-2} -E \|P_n-P\|_{\mathcal G_k} \ge c_0 2^m$$ for some fixed $c_0$ precisely when $2^m$ is of larger magnitude than $(2^m)^{\frac{1}{2}-\frac{1}{4s}} n^{-1/2} + (2^m)^{-1/2s}n^{-1}$, equivalent to $2^m \ge c'' n^{-2s/(2s+1)}$ for some $c''>0$, which is satisfied since $2^m \ge c' \rho_n^2 \ge c'' n^{-2s/(2s+1)}$ if $L_0$ is large enough, by hypothesis on $\rho_n$. We can thus rewrite the last probability in (\ref{sup}) as 
$$ \sum_{m \in \mathbb Z : c'\rho_n^2 \le 2^m \le C} {\Pr}_f \left\{ n\|P_n-P\|_{\mathcal G_{m, J_n}}- nE\|P_n-P\|_{\mathcal G_{m, J_n}} > c_0n2^m \right\}.$$
To this expression we can apply Talagrand's inequality (\ref{bous}), noting that the supremum over $\mathcal G_{m, J_n}$ can be realised, by continuity, as one over a countable subset of $\Sigma$, and since $\Sigma$ is uniformly bounded by $\sup_{f \in \Sigma(t, B)}\|f\|_\infty \le U \equiv U(t,B)$. Renormalising by $U$ and using (\ref{bous}), (\ref{weakvar}), (\ref{mom}) we can bound the expression in the last display, up to multiplicative constants, by
\begin{eqnarray*}
\sum_{m \in \mathbb Z: c'\rho_n^2 \le 2^m \le C} \exp \left\{-c_1 \frac{n^2(2^m)^2}{n2^m + nE \|P_n-P\|_{\mathcal G_{m, J_n}}+ n2^m } \right\} &\le& \sum_{m \in \mathbb Z: c'\rho_n^2 \le 2^m \le C} e^{-c_2 n2^m} \\
&\le & c_3 e^{-c_4 n \rho_n^2}
\end{eqnarray*}
since $2^m \ge c'\rho_n^2 >> n^{-1}$, which completes the proof.
\end{proof}

\subsection{Proof of Theorem \ref{adapt}}

\begin{proof} We construct a standard Lepski type estimator: choose integers $j_{\min}, j_{\max}$ such that $J_0 \le j_{\min} < j_{\max}$, $$2^{j_{\min}} \simeq n^{1/(2R+1)} ~~\textrm{and} ~~ 2^{j_{\max}} \simeq n^{1/(2r+1)}$$ and define the grid $$ \mathcal J := \mathcal J_n = [j_{\min}, j_{\max}] \cap \mathbb N.$$ Let $f_n(j) \equiv f_n(j,\cdot)=\int_0^1 K_j(\cdot,y)dP_n(y)$ be a linear wavelet estimator based on wavelets of regularity $S>R$. To simplify the exposition we prove the result for $\|f\|_\infty$ known, otherwise the result follows from the same proof, with $\|f\|_\infty$ replaced by $\|f_n(j_{\max})\|_\infty$, a consistent estimator for $\|f\|_\infty$ that satisfies sufficiently tight uniform exponential error bounds  (using inequality (26) in \cite{GN11} and proceeding as in Step (II) on p.1157 in \cite{GN10b}). Set
\begin{equation}
\bar j_n = \min \bigg \{ j \in \mathcal J: \|f_n(j) - f_n(l)\|_2^2 \le C(S) (\|f\|_\infty \vee 1) \frac{2^l}{n} ~~ \forall l>j, l\in {\mathcal J} \bigg \}
\label{htf2}
\end{equation}
where $C(S)$ is a large enough constant, to be chosen below, in dependence of the wavelet basis. The adaptive estimator is $\hat f_n = f_n(\bar j_n)$. We shall need the standard estimates
\begin{equation}
E\|f_n(j) - Ef_n(j)\|_2^2 \leq D \frac{2^j }{n} := D \sigma^2 (j,n)
\label{var}
\end{equation}
and, for $f \in W^s, s \in [r,R]$,
\begin{equation}
\|E f_n(j) - f\|_2 \leq 2^{-js}  D' \|f\|_{s,2} := B(j, f)
\label{bias}
\end{equation}
for constants $D, D'$ that depend only on the wavelet basis and on $r,R$. Define $j^*:=j^*(f)$ by
\begin{equation*} 
j^*=\min\left\{j\in {\cal J}: B(j,f) \le \sqrt D \sigma(j,n) \right\}
\end{equation*} 
so that, for every $f \in \Sigma(s,B)$ and $D''=D''(D,D')$
\begin{equation} \label{bal}
D^{-1} B^2(j^*,f) \le \sigma^2(j^*,n) \le D'' \|f\|_{s,2}^{2/(2s+1)} n^{-2s/(2s+1)} \le D'' B^{2/(2s+1)} n^{-2s/(2s+1)}.
\end{equation}
We will consider the cases $\{\bar j_n \leq j^* \}$ and $\{\bar j_n > j^* \}$ separately. First, by the definition of $\bar j_n, j^*$ and (\ref{var}), (\ref{bias}), (\ref{bal}),
\begin{eqnarray} \label{lowvar}
E \left\|f_n(\bar j_n) - f \right\|^2_2 I_{\{\bar j_n \le j^* \}} &=& E \left( \|f_n (\bar j_n) - f_n (j^*) \|^2_2 + E\|f_n (j^*) - f\|^2_2 \right) I_{\{\bar j_n \le j^* \}} \notag \\
& \le & C(S)(\|f\|_\infty \vee 1)\frac{2^{j^*}}{n} + C' \sigma^2(j^*,n) \le C'' B^{2/(2s+1)} n^{-2s/2s+1} \notag
\end{eqnarray}
for $C''=C''(D,D', S,U)$, which is the desired bound. On the event $\{\bar j_n > j^* \}$ we have, using (\ref{var}) and the definition of $j^*$,
\begin{eqnarray*}
E \left\|f_n(\bar j_n) - f\right\|_2 I_{\{\hat j_n > j^* \}}  & \leq & \sum_{j \in \mathcal{J}:j > j^*} \left(E \left\|f_n(j) -f\right\|_2^2 \right)^{1/2} ~  \left(EI_{\{\hat j_n = j\}}\right)^{1/2}  \\
& \leq &  \sum_{j \in \mathcal{J}: j > j^*} C''' \sigma(j,n) \cdot \sqrt{ {\Pr}_f\{\hat j_n = j\}} \\
&\le & C'''' \sum_{j \in \mathcal J: j>j^*} \sqrt{ {\Pr}_f\{\hat j_n = j\}}
\end{eqnarray*}
since $\sup_{j \in \mathcal J}\sigma(j,n) = \sigma(j_{\max},n)$ is bounded in $n$. Now pick any $j \in \mathcal{J}$ so that $j > j^*$ and denote by $j^-$ the previous element in the grid (i.e. $j^-= j-1$). One has, by definition of $\bar j_n$,
\begin{equation} \label{pr}
{\Pr}_f \{\bar j_n= j\}  \leq  \sum_{l \in \mathcal{J}: l \ge j} {\Pr}_f \left \{ \left\|f_n(j^-) - f_n(l) \right\|_2 > \sqrt {C(S) (\|f\|_\infty \vee 1) \frac{2^l}{n}}\right\},
\end{equation}
and we observe that, by the triangle inequality,
\begin{equation*}
\left\|f_n(j^-) - f_n(l) \right \|_2 \leq \left\|f_n(j^-) - f_n(l) - Ef_n(j^-) + Ef_n(l)  \right\|_2 + B(j^-, f) + B(l, f), 
\end{equation*} 
where, $$ B(j^-, f) + B(l, f) \leq 2B(j^*, f) \leq c \sigma(j^*,n) \leq c' \sigma(l,n) $$ by definition of $j^*$ and since $l>j^- \geq j^*$.
Consequently, the probability in (\ref{pr}) is bounded by
\begin{equation} \label{c1} 
{\Pr}_f \left \{ \left\|f_n(j^-) -  f_n(l) -Ef_n(j^-) + Ef_n(l) \right\|_2 > (\sqrt{C(S)(\|f\|_\infty \vee 1)}-c') \sigma(l,n) \right \},
\end{equation}
and by inequality (\ref{tal2}) above this probability is bounded by a constant multiple of $e^{-d2^l}$ if we choose $C(S)$ large enough. This gives the overall bound $$\sum_{l \in \mathcal J: l \geq j} c''e^{-d2^l} \le d'e^{-d''2^{j_{\min}}},$$ which is smaller than a constant multiple times $B^{1/(2s+1)} n^{-s/(2s+1)}$, uniformly in $s \in [r,R], n \in \mathbb N$ and for $B \ge 1$, by definition of $j_{\min}$. This completes the proof.
\end{proof}

\subsection{Proof of Theorem \ref{rv}}

\begin{proof}

A) Suppose for simplicity that the sample size is $2n$, and split the sample into two halves with index sets $\mathcal S^1, \mathcal S^2$, of equal size $n$, write $E_1, E_2$ for the corresponding expectations, and $E=E_1E_2$. Let $\hat f_n= f_n(\bar j_n)$ be the adaptive estimator from the proof of Theorem \ref{adapt} based on the sample $\mathcal S^1$. One shows by a standard bias-variance decomposition, using $\bar j_n \in \mathcal J$ and $\|K_j(f)\|_{r,2} \le \|f\|_{r,2}$, that for every $\varepsilon>0$ there exists a finite positive constant $B'=B'(\varepsilon, B_0)$ satisfying $$\inf_{f \in \Sigma (r,B_0)}{\Pr}_f\{\|\hat f_n\|_{r,2} \le B'\} \ge 1 - \varepsilon.$$ It therefore suffices to prove the theorem on the event $\{\|\hat f_n\|_{r,2} \le B'\}$. For a wavelet basis of regularity $S>R$ and for $J_n \ge J_0$ a sequence of integers such that $2^{J_n} \simeq n^{1/(2r+1/2)}$, define the $U$-statistic
\begin{equation}\label{ustat0}
U_n(\hat f_n)=\frac{2}{n(n-1)} \sum_{i<j, i,j \in \mathcal S^2} \sum_{l=J_0}^{J_n-1} \sum_{k \in \mathcal Z_l} (\psi_{lk}(X_i)-\langle \psi_{lk}, \hat f_n\rangle)(\psi_{lk}(X_j)-\langle \psi_{lk}, \hat f_n \rangle)
\end{equation}
which has expectation $$E_2 U_n(\hat f_n) = \sum_{l=J_0}^{J_n-1} \sum_{k \in \mathcal Z_l}\langle \psi_{lk}, f- \hat f_n \rangle^2 = \|\Pi_{V_{J_n}}(f-\hat f_n)\|_2^2.$$ Using Chebychev's inequality and that, by definition of the norm (\ref{sobolev}) $$\sup_{h \in \Sigma(r,b)}\|\Pi_{V_{J_n}}(h)-h\|_2^2 \le c(b) 2^{-2J_nr}$$ for every $0<b<\infty$ and some finite constant $c(b)$, we deduce
\begin{eqnarray*} \label{cov}
&&\inf_{f \in \Sigma(r,B_0)} {\Pr}_{f,2} \left\{U_n(\hat f_n) - \|f-\hat f_n\|_2^2  \ge -(c(B_0)+c(B'))2^{-2J_nr} -z(\alpha) \tau_n(f) \right\} \\
&& \ge \inf_{f \in \Sigma(r,B_0)} {\Pr}_{f,2} \left\{U_n(\hat f_n) - \|\Pi_{V_{J_n}}(f-\hat f_n)\|_2^2  \ge  -z(\alpha) \tau_n(f) \right\}  \\
&& \ge 1 - \sup_{f \in \Sigma(r,B_0)}\frac{Var_2(U_n(\hat f_n)-E_2U_n(\hat f_n))}{(z(\alpha) \tau_n(f))^2}.
\end{eqnarray*}
We now show that the last quantity is greater than or equal to $1-z(\alpha)^{-2} \ge 1- \alpha$ for quantile constants $z(\alpha)$ and with $$\tau^2_n(f)= \frac{C(S)2^{J_n}\|f\|^2_\infty}{n(n-1)}+\frac{4\|f\|_\infty}{n}\|\Pi_{V_{J_n}}(f-\hat f_n)\|_2^2,$$ which in turn gives the honest confidence set under $\Pr$
\begin{equation} \label{conf0} 
C_n(\|f\|_\infty, B_0) = \left \{f: \|f-\hat f_n\|_2 \le \sqrt{z_\alpha \tau_n(f) + U_n(\hat f_n) + (c(B_0)+c(B'))2^{-2{J_n}r}}  \right\}.
\end{equation}
We shall comment on the role of the constants $\|f\|_\infty, c(B_0), C(B')$ at the end of the proof, and establish the last claim first: note that the Hoeffding decomposition for the centered $U$-statistic with kernel $$R(x,y)= \sum_{l=J_0}^{J_n-1} \sum_{k \in \mathcal Z_l}  (\psi_{lk}(x)-\langle \psi_{lk}, \hat f_n\rangle)(\psi_{lk}(y)-\langle \psi_{lk}, \hat f_n \rangle)$$ is (cf.~the proof of Theorem 4.1 in \cite{RV06}) 
$$U_n(\hat f_n)-E_2U_n(\hat f_n) = \frac{2}{n} \sum_{i=1}^n (\pi_1R)(X_i) + \frac{2}{n(n-1)}\sum_{i<j}(\pi_2R)(X_i, X_j) \equiv L_n + D_n$$
where 
\begin{equation*}
(\pi_1R)(x) =\sum_{l=J_0}^{J_n-1} \sum_{k \in \mathcal Z_l} (\psi_{lk}(x)-\langle \psi_{lk}, f \rangle) \langle \psi_{lk}, f-\hat f_n \rangle
\end{equation*}
and $$(\pi_2R)(x,y)= \sum_{l=J_0}^{J_n-1} \sum_{k \in \mathcal Z_l} (\psi_{lk}(x)-\langle \psi_{lk}, f \rangle)(\psi_{lk}(y)-\langle \psi_{lk}, f \rangle)$$
The variance of $U_n(\hat f_n)-E_2U_n(\hat f_n)$ is the sum of the variances of the two terms in the Hoeffding decomposition. For the linear term we bound the variance $Var_2(L_n)$ by the second moment, using orthonormality of the $\psi_{lk}$s,
\begin{equation*}
\frac{4}{n} \int \left(\sum_{l=J_0}^{J_n-1} \sum_{k \in \mathcal Z_l} \psi_{lk}(x) \langle \psi_{lk}, \hat f_n -f \rangle \right)^2 f(x) dx \le \frac{4 \|f\|_\infty}{n} \sum_{l=J_0}^{J_n-1} \sum_{k \in \mathcal Z_l} \langle  \psi_{lk}, \hat f_n-f \rangle^2,
\end{equation*}
which equals the second term in the definition of $\tau^2_n(f)$. For the degenerate term we can bound $Var_2(D_n)$ analogously by the second moment of the uncentered kernel (cf.~after (\ref{bd0})), i.e., by
\begin{equation*}
\frac{2}{n(n-1)} \int \left(\sum_{l=J_0}^{J_n-1} \sum_{k \in \mathcal Z_l} \psi_{lk}(x) \psi_{lk}(y) \right)^2 f(x) dx f(y)dy \le \frac{C(S) 2^{J_n} \|f\|^2_\infty}{n(n-1)}, 
\end{equation*}
using orthonormality and the cardinality properties of $\mathcal Z_l$. 

The so constructed confidence set has an adaptive expected maximal diameter: let $f \in \Sigma(s,B)$ for some $s \in [r,R]$ and some $1 \le B \le B_0$. The nonrandom terms are of order $$\sqrt{c(B_0)+c(B')}2^{-J_nr} + \|f\|_\infty^{1/2}2^{J_n/4}n^{-1/2} \le C(S, B_0, B', r, U) n^{-r/(2r+1/2)}$$ which is $o(n^{-s/(2s+1)})$ since $s \le R < 2r$. The random component of $\tau_n(f)$ has order $\|f\|_\infty^{1/4} n^{-1/4}E_1\|\Pi_{V_{J_n}}(\hat f_n - f)\|_2^{1/2}$ which is also $o(n^{-s/(2s+1)})$ for $s<2r$, since $\Pi_{V_{J_n}}$ is a projection operator and since $\hat f_n$ is adaptive, as established in Theorem \ref{adapt}. Moreover, by Theorem \ref{adapt} and again the projection properties, $$EU_n(\hat f_n) = E_1\|\Pi_{V_{J_n}}(\hat f_n-f)\|_2^2 \le E_1\|\hat f_n -f\|_2^2 \le c B^{2/(2s+1)}n^{-2s/(2s+1)}.$$ The term in the last display is the leading term in our bound for the diameter of the confidence set, and shows that $C_n$ adapts to both $B$ and $s$ in the sense of Definition \ref{cc}, using Markov's inequality.

The confidence set $C_n(\|f\|_\infty, B_0)$ is not feasible if $B_0$ and $\|f\|_\infty$ are unknown, so in particular under the assumptions of Theorem \ref{rv}, but $C_n$ independent of $B_0, \|f\|_\infty$ can be constructed as follows: we replace $c(B_0)+c(B')$ in the definition of (\ref{conf0}) by a divergent sequence of positive real numbers $c_n$, which can still be accommodated in the diameter estimate from the last paragraph since $n^{-2r/(2r+1/2)}c_n$ is still $o(n^{-2s/(2s+1)})$ as long as $s \le R<2r$ for $c_n$ diverging slowly enough (e.g., like $\log n$). Define thus the confidence set
\begin{equation} \label{conf}
C_n = \left \{f: \|f-\hat f_n\|_2 \le \sqrt{z_\alpha \tau_n(f) + U_n(\hat f_n) + c_n2^{-2Jr}}  \right\},
\end{equation}
with $\|f\|_\infty$ replaced by $\|f_n(j_{\max})\|_\infty$ in all expressions where $\|f\|_\infty$ occurs. As stated before (\ref{htf2}), $\|f_n(j_{\max})\|_\infty$ concentrates around $\|f\|_\infty$ with exponential error bounds, so that the sufficiency part of Theorem \ref{rv} then holds for this $C_n$ with slightly increased $z_\alpha$. 

\medskip

B) Necessity of $R \le 2r$ follows immediately from Part B of Theorem \ref{fix}. That $R<2r$ is also necessary is proved in Subsection \ref{lpt} below.
\end{proof}

\subsection{Proof of Theorem \ref{fix}}

\begin{proof}
That an $L^2$-adaptive confidence set exists when $s \le 2r$ follows from Theorem \ref{rv}; The case $s<2r$ is immediate, and the case $s=2r$ follows using the confidence set (\ref{conf0}). This set is feasible since, under the hypotheses of Theorem \ref{fix}, $B=B_0$ is known, as is $B'$ and the upper bound for $\|f\|_\infty$ (cf.~(\ref{supbd})). It is further adaptive since $n^{-r/(2r+1/2)}=n^{-s/(2s+1)}$ for $s=2r$. 

For part Aii we use the test $\Psi_n$ from Proposition \ref{test} with $\Sigma=\Sigma(s), t=r,$ and define a confidence ball as follows. Take $\hat f_n=f_n(\bar j_n)$ to be the adaptive estimator from the proof of Theorem \ref{adapt}, and let, for $0<L'<\infty$,
\begin{eqnarray*}
C_n=\begin{cases}\{f \in \Sigma(r): \|f-\hat f_n\|_2 \le L'n^{-s/(2s+1)}\}&\text{if} ~\Psi_n=0\\
\{f \in \Sigma(r): \|f-\hat f_n\|_2 \le L'n^{-r/(2r+1)}\}&\text{if}~\Psi_n=1 \end{cases}
\end{eqnarray*}
We first prove that $C_n$ is honest for $\Sigma(s) \cup \tilde \Sigma(r, \rho_n)$ if we choose $L'$ large enough. For $f \in \Sigma(s)$ we have from Theorem \ref{adapt}, by Markov's inequality,
\begin{eqnarray*}
\inf_{f \in \Sigma(s)}{\Pr}_f \left\{f \in C_n \right\}&\ge& 1- \sup_{f \in \Sigma(s)}{\Pr}_f \left \{ \|\hat f_n-f\|_2 >  L'n^{-s/(2s+1)}  \right\} \\
& \ge & 1-\frac{n^{s/(2s+1)}}{L'} \sup_{f \in \Sigma(s)}E_f \|\hat f_n -f\|_2 \\
&\ge& 1-\frac{c(B,s,r)}{L'}
\end{eqnarray*}
which can be made greater than $1-\alpha$ for any $\alpha>0$ by choosing $L'$ large enough depending only on $B, \alpha, r,s$. When $f \in \tilde \Sigma(r, \rho_n)$, using again Markov's inequality
\begin{equation*}
\inf_{f \in \tilde \Sigma(r, \rho_n)}{\Pr}_f \left\{f \in C_n\right\} \ge 1 - \frac{\sup_{f \in \Sigma(r)}E_f\|\hat f_n -f\|_2}{L'n^{-r/(2r+1)}} - \sup_{f \in \tilde \Sigma(r, \rho_n)}{\Pr}_f\{\Psi_n =0\}.
\end{equation*}
The first subtracted term can be made smaller than $\alpha/2$ for $L'$ large enough as before. The second subtracted term can also be made less than $\alpha/2$ using Proposition \ref{test} and the remark preceding it, choosing $M$ and $d_n$ to be large but also bounded in $n$. This proves that $C_n$ is honest. We now turn to adaptivity of $C_n$: by the definition of $C_n$ we always have $|C_n| \le L'n^{-r/(2r+1)}$, so the case $f \in \tilde \Sigma(r, \rho_n)$ is proved. If $f \in \Sigma(s)$ then using Proposition \ref{test} again, for $M, d_n$ large enough depending on $\alpha'$ but bounded in $n$, $${\Pr}_f\{|C_n| > L'n^{-s/(2s+1)}\} = {\Pr}_f\{\Psi_n =1\} \le \alpha',$$ which completes the proof of part A.

To prove part B of Theorem \ref{fix} we argue by contradiction and assume that the limit inferior equals zero. We then pass to a subsequence of $n$ for which the limit is zero, and still denote this subsequence by $n$. Let $f_0\equiv 1 \in \Sigma(s)$, suppose $C_n$ is adaptive and honest for $\Sigma(s) \cup \tilde\Sigma(r, \rho_n)$ for every $\alpha, \alpha'$, and consider testing $$H_0: f=f_0 ~~~\text{against}~~~H_1: f \in \tilde \Sigma(r, \rho_n)$$ where $\rho_n =o(n^{-r/(2r+1/2)})$. Since $s>2r$ we may assume $n^{-s/(2s+1)}=o(\rho_n)$ (otherwise replace $\rho_n$ by $\rho_n' \ge \rho_n$ s.t. $n^{-s/(2s+1)}=o(\rho_n')$). Accept $H_0$ if $C_n \cap \tilde \Sigma(r, \rho_n)$ is empty and reject otherwise, formally $$\Psi_n = 1\{C_n \cap \tilde \Sigma(r, \rho_n) \neq \emptyset\}.$$ The type-one errors of this test satisfy
\begin{eqnarray*}
E_{f_0}\Psi_n &=& {\Pr}_{f_0}\left\{C_n \cap \tilde \Sigma(r, \rho_n) \neq \emptyset \right\} \\
&\le & {\Pr}_{f_0} \{f_0 \in C_n, |C_n| \ge \rho_n\} + {\Pr}_{f_0} \{f_0 \notin C_n\} \\
&\le & \alpha + \alpha' + r_n \to \alpha + \alpha'
\end{eqnarray*}
as $n \to \infty$ by the hypothesis of coverage and adaptivity of $C_n$. The type-two errors satisfy, by coverage of $C_n$, as $n \to \infty$
$$E_f(1-\Psi_n) = {\Pr}_f \{C_n \cap \tilde \Sigma(r, \rho_n)=\emptyset\} \le {\Pr}_{f} \{f \notin C_n\} \le \alpha + r_n \to \alpha,$$ uniformly in $f \in \tilde \Sigma(r, \rho_n)$. We conclude that this test satisfies $$\limsup_n\left[E_{f_0}\Psi_n + \sup_{f \in H_1}E_f(1-\Psi_n)\right] \le 2\alpha + \alpha'$$ for arbitrary $\alpha, \alpha'>0$. For $\alpha, \alpha'$ small enough this contradicts (the proof of) Theorem 1i in \cite{I86}, which implies that the limit inferior of the term in brackets in the last display, even with an infimum over all tests, exceeds a fixed positive constant. Indeed, the alternatives (6) in \cite{I86} can be taken to be $$f_i(x) = 1 + \epsilon 2^{-j_n (r+1/2)} \sum_{k \in \mathcal Z_{j_n}} \beta_{ik} \psi_{j_n k}(x), ~~~~~i=1, \dots, 2^{2^{j_n}},$$ for $\epsilon>0$ a small constant, $\beta_{ik} = \pm 1$, and with $j_n$ such that $2^{j_n} \simeq n^{1/(2r+1/2)}$. Since $$\inf_{g \in \Sigma(s)}\|f_i-g\|_2 \ge \sqrt{\sum_{l\ge j_n, k}\langle f_i, \psi_{lk}\rangle^2} - \sup_{g \in \Sigma(s)} \sqrt{\sum_{l\ge j_n, k}\langle g, \psi_{lk}\rangle^2} \ge c\epsilon n^{-r/(2r+1/2)}$$ for every $\epsilon>0$, some $c>0$ and $n$ large enough, these alternatives are also contained in our $H_1$, so that the proof of the lower bound Theorem 1i in \cite{I86} applies also in the present situation.
\end{proof}

\subsection{Proof of Theorem \ref{cont}}

We shall write $\Sigma(s)$ for $\Sigma(s,B_0)$ and $\tilde \Sigma_n(s)$ for $\tilde \Sigma(s, \rho_n(s))$ in this proof, and we write $\tilde \Sigma_n(s_N)$ also for $\Sigma(s_N)$ in slight abuse of notation. For $i=1,\dots, N,$ let $\Psi(i)$ be the test from (\ref{stat}) with $\Sigma=\Sigma(s_{i+1})$ and $t=s_i$. Starting from the largest model we first test $H_0: f \in \Sigma(s_2)$ against $H_1: f \in \tilde \Sigma_n(s_1)$, accepting $H_0$ if $\Psi(1)=0$. If $H_0$ is rejected we set $\hat s_n = s_1=r$, otherwise we proceed to test $H_0: f \in \Sigma (s_3)$ against $H_1: f \in \tilde \Sigma_n(s_2)$ using $\Psi(2)$ and iterating this procedure downwards we define $\hat s_n$ to be the first element $s_i$ in $\mathcal S$ for which $\Psi(i)=1$ rejects. If no rejection occurs we set $\hat s_n$ equal to $s_N$, the last element in the grid. 

For $f \in \mathcal P_n(M, \mathcal S)$ define the unique $s_{i_0}:=s_{i_0}(f)= \{s \in \mathcal S: f \in \tilde \Sigma_n(s) \}$. We now show that for $M$ large enough
\begin{equation} \label{consist}
\sup_{f \in \mathcal P_n(M, \mathcal S)}{\Pr}_f \{\hat s_n \ne s_{i_0}(f)\} < \max(\alpha, \alpha')/2.
\end{equation}
Indeed, if $\hat s_n < s_{i_0}$ then the test $\Psi(i)$ has rejected for some $i < i_0$. In this case $f \in \tilde \Sigma_n(s_{i_0}) \subset \Sigma(s_{i_0}) \subseteq \Sigma(s_{i+1})$ for every $i<i_0$, and thus,
\begin{eqnarray*}
{\Pr}_f\{\hat s_n< s_{i_0}\} &=& {\Pr}_f\left\{\bigcup_{i < i_0} \{\Psi(i)=1\} \right\} \le \sum_{i<i_0} \sup_{f \in \Sigma(s_{i+1})}E_f\Psi(i) \\
&\le&  C(N) e^{-cd_n^2} < \max(\alpha, \alpha')/2
\end{eqnarray*}
 using Proposition \ref{test} and the remark preceding it, choosing $M$ and $d_n$ to be large but also bounded in $n$. On the other hand if $\hat s_n > s_{i_0}$ (ignoring the trivial case $s_{i_0} =s_N$) then $\Psi(i_0)$ has accepted despite $f \in \tilde \Sigma_n(s_{i_0})$. Thus $${\Pr}_f\{\hat s_n > s_{i_0}\} \le \sup_{f \in \tilde \Sigma_n(s_{i_0})}E_f (1-\Psi(i_0)) \le C e^{-cd_n^2} \le \max(\alpha, \alpha')/2 $$ again by Proposition \ref{test}, for $M, d_n$ large enough.

Denote now by $C_n(s_i)$ the confidence set (\ref{conf0}) constructed in the proof of Theorem \ref{rv} with $r$ there being $s_i$, with $R=2s_i=s_{i+1}$, with $\|f\|_\infty$ replaced by $U$ and with $z_\alpha$ such that the asymptotic coverage level is $\alpha/2$ for any $f \in \Sigma(s_i)$. We then set $C_n=C_n(\hat s_n)$, which is a feasible confidence set as $B_0, r, U$ are known under the hypotheses of the theorem. We then have, from the proof of Theorem \ref{rv}, uniformly in $f \in \tilde \Sigma_n(s_{i_0}) \subset \Sigma(s_{i_0})$, $${\Pr}_f\{f \in C_n(\hat s_n)\} \ge {\Pr}_f\{f \in C_n(s_{i_0})\} - \alpha/2 \ge 1- \alpha.$$ Moreover, if $f \in \Sigma (s, B) \cap \tilde \Sigma_n(s_{i_0})$ for some $1 \le B \le B_0$ and for either $s \in [s_{i_0}, s_{i_0+1})$ or $s \in [s_N, R]$ (in case $s_{i_0}=s_N$), the expected diameter of $C_n$ satisfies, by the estimates in the proof of Theorem \ref{rv},
\begin{eqnarray*}
&& {\Pr}_f\{|C_n(\hat s_n)| > C B^{2/(2s+1)}n^{-s/(2s+1)}\} \\
&& \le {\Pr}_f\{|C_n(s_{i_0})| > C B^{2/(2s+1)}n^{-s/(2s+1)}\} + \alpha'/2 \\
&& \le \alpha'
\end{eqnarray*}
for $C$ large enough, so that this confidence set is adaptive as well, which completes the proof.

\subsection{Proof of Theorem \ref{imp}}

\begin{proof}

Suppose such $C_n$ exists. We will construct functions $f_m \in
W^s, m = 0, 1, \dots,$ and a further function $f_\infty \in W^r$, which
serve as hypotheses for $f$. For each $m \in \mathbb N$, we
will ensure that, at some time $n_m$, $C_{n_m}$ cannot distinguish between
$f_m$ and $f_\infty$, and is too small to contain both simultaneously. We
will thereby obtain a subsequence $n_m$ on which, for  \(\delta = \tfrac15(1 - 
  2\alpha),\)  \[\sup_m \Pr_{f_\infty} \{f_\infty \in C_{n_m}\} \le 1 - \alpha - \delta,\]
  contradicting our assumptions on \(C_n.\)

For \(m = 0, 1, 2, \dots, \infty,\) construct functions $f_0=1$, 
  \[f_m = 1 + \varepsilon \sum_{i=1}^m \sum_{k \in \mathcal Z_{j_i}} 
  2^{-j_i(r+1/2)} \beta_{ik} \psi_{j_ik}.\]
  where $\varepsilon>0$ is a constant, and the parameters \(j_1, j_2, \ldots \in \mathbb N\),  \(\beta_{ik} = \pm 1\) are chosen inductively  satisfying \(j_i/j_{i-1} \ge 
  1 + 1/2r\).  Pick \(\varepsilon > 0\) small 
  enough that \(\|f_m - f_{m-1}\|_\infty \le 2^{-(m+1)}\) for all \(m < 
  \infty,\) and any choice of \(j_i, \beta_{ik}.\) Then
  \[f_m = 1 + \sum_{i=1}^m (f_i - f_{i-1}) \ge \tfrac12,\]
  and \(\int f_m = \langle 1, f_m \rangle = 1,\) so the \(f_m\) are densities.  
 By (\ref{sobolev}), \(f_m \in W^r,\) and for \(m < \infty,\) also \(f_m \in W^s.\) 
 
  We have already defined $f_0$; for convenience let $n_0 = 1$. Inductively, suppose we have defined \(f_{m-1}, n_{m-1}.\)  For \(n_m >   n_{m-1}\) and \(D>0 \) large enough depending only on $f_{m-1}$, we have:
  \begin{enumerate}
    \item \(\Pr_{f_{m-1}}\{f_{m-1} \not\in C_{n_m}\} \le \alpha + \delta\); and
    \item \(\Pr_{f_{m-1}}\{|C_{n_m}| \ge Dr_{n_m}\} \le \delta.\)
  \end{enumerate}
  Setting $$T_n = 1(\exists\ f \in C_n, \|f - f_{m-1}\|_2 \ge 2Dr_n),$$ we 
  then have
  \begin{equation}
    \label{eq:c-accurate}
    \Pr_{f_{m-1}}\{T_{n_m}=1\} \le \Pr_{f_{m-1}}\{f_{m-1} \not \in C_{n_m}\} + 
    \Pr_{f_{m-1}}\{|C_{n_m}| \ge Dr_{n_m}\} \le \alpha + 2\delta.
  \end{equation}
  We claim it is possible to choose $j_m, \beta_{mk}$ and $n_m$, depending only
on $f_{m-1}$ so that also: 1. if $m>1$,
  \begin{equation}
    \label{eq:f-separate}
    3Dr_{n_m} \le \|f_m - f_{m-1}\|_2 \le \tfrac14 \|f_{m-1} - 
    f_{m-2}\|_2,
  \end{equation}
  and 2. for any further choice of \(j_i, \beta_{ik},\)
  \begin{equation}
    \label{eq:f-indistinguishable}
    \Pr_{f_\infty}\{T_{n_m} = 0\} \ge 1 - \alpha - 4\delta.
  \end{equation}
  We may then conclude that, since all further choices will satisfy 
  \eqref{eq:f-separate},
  \[\|f_\infty - f_{m-1}\|_2 \ge \|f_m - f_{m-1}\|_2 - 
  \sum_{i=m+1}^\infty \|f_i - f_{i-1}\|_2 \ge 2Dr_{n_m},\]
  so
  \[\Pr_{f_\infty}\{f_\infty \in C_{n_m}\} \le \Pr_{f_\infty}\{T_{n_m} = 1\} \le 
  \alpha + 4\delta = 1 - \alpha - \delta\]
  as required.

  It remains to verify the claim. For \(j \ge (1 + 1/2r)j_{m-1},\) \(\beta_k = 
  \pm 1,\) set
  \[g_\beta = \varepsilon 2^{-j(r+1/2)} \sum_{k \in \mathcal Z_j} \beta_k 
  \psi_{jk},\]
  and \(f_\beta = f_{m-1} + g_\beta.\) Allowing \(j \to \infty,\) set \[n \sim 
  C2^{j(2r+1/2)},\] for \(C > 0\) to be determined. Then $$\|g_\beta\|_2 = \varepsilon 2^{-jr} \approx n^{-r/(2r+1/2)},$$ so for \(j\) large enough, \(f_\beta\) satisfies \eqref{eq:f-separate} with any 
  choice of \(\beta.\)

  The density of \(X_1, \dots, X_n\) under \(f_\beta,\) w.r.t.\ under \(f_{m-1},\) 
  is
  $$Z_\beta = \prod_{i=1}^n \frac{f_\beta}{f_{m-1}}(X_i).$$
  Set \(Z = 2^{-j}\sum_\beta Z_\beta,\) so \(E_{f_{m-1}}[Z] = 1,\) and
  \begin{align*}
    E_{f_{m-1}}[Z^2] &= 2^{-2j} \sum_{\beta, \beta'} \prod_{i=1}^n 
    E_{f_{m-1}}\left[ \frac{f_\beta f_{\beta'}}{f_{m-1}^2}(X_i)\right]\\
    &= 2^{-2j} \sum_{\beta, \beta'} \left \langle 
    \frac{f_\beta}{\sqrt{f_{m-1}}}, \frac{f_{\beta'}}{\sqrt{f_{m-1}}} \right 
    \rangle^n\\
    &= 2^{-2j} \sum_{\beta, \beta'} \left(1 + \left \langle 
    \frac{g_\beta}{\sqrt{f_{m-1}}}, \frac{g_{\beta'}}{\sqrt{f_{m-1}}} \right 
    \rangle\right)^n\\
    &\le 2^{-2j} \sum_{\beta, \beta'} (1 + 2 \langle \beta, \beta' 
    \rangle)^n\\
    &= E[(1 + \varepsilon^22^{1-j(2r+1)} Y)^n],\\
    \intertext{where \(Y = \sum_{i=1}^{2^j} R_i,\) for i.i.d.\ Rademacher random 
    variables \(R_i,\)}
    &\le E[\exp(n\varepsilon^2 2^{1-j(2r+1)} Y)]\\
    &= \cosh\left(D2^{-j/2}(1 + o(1))\right)^{2^j},\\
    \shortintertext{as \(j \to \infty,\) for some \(D > 0,\)}
    &= \left(1 + D^2 2^{-j} (1 + o(1))\right)^{2^j}\\
    &\le \exp\left(D^2(1 + o(1))\right)\\
    &\le 1 + \delta^2,
  \end{align*}
  for \(j\) large, \(C\) small. Hence
  \(E_{f_{m-1}}[(Z - 1)^2] \le \delta^2,\)
  and we obtain
  \begin{align*}
    \Pr_{f_{m-1}}\{T_n = 1\} + \max_\beta \Pr_{f_\beta}\{T_n = 0\}
    &\ge \Pr_{f_{m-1}}\{T_n=1\} + 2^{-j}\sum_\beta \Pr_{f_\beta}\{T_n = 0\}\\
    &= 1 + E_{f_{m-1}}[(Z-1)1(T_n = 0)]\\
    &\ge 1 - \delta.
  \end{align*}

  Set \(f_m = f_\beta,\) for \(\beta\) maximizing this expression. The density of 
  \(X_1, \dots, X_n\) under \(f_\infty,\) w.r.t.\ under \(f_m,\) is
  \[Z' = \prod_{i=1}^n \frac{f_\infty}{f_m}(X_i).\]
  Now, \(E_{f_m}[Z'] = 1,\) and
  \[\norm{f_\infty - f_m}_2^2 = \sum_{i=m+1}^\infty \varepsilon^2 2^{-2j_ir} 
  \le E'2^{-2j_{m+1}r} \le E'2^{-j(2r+1)},\]
  for some constant \(E' > 0,\) so similarly
  \begin{align*}
    E_{f_m}[{Z'}^2] &\le (1 + 2\norm{f_\infty - f_m}_2^2)^n\\
    &\le (1 + E' 2^{1 - j(2r+1)})^n\\
    &\le \exp(E'n2^{1 - j(2r+1)})\\
    &= \exp\left(F2^{-j/2}(1+o(1))\right),\\
    \shortintertext{for some \(F > 0,\)}
    &\le 1 + \delta^2,
  \end{align*}
  for \(j\) large. Hence
  \(E_{f_{m}}[(Z'-1)^2] \le \delta^2,\)
  and
  \begin{align*}
    \Pr_{f_{m-1}}\{T_n = 1\} + \Pr_{f_\infty}\{T_n = 0\}
    &= \Pr_{f_{m-1}}\{T_n=1\} + E_{f_m}[Z'1(T_n=0)]
    \\&\ge 1 - \delta + E_{f_m}[(Z'-1)1(T_n = 0)]
    \\&\ge 1 - 2\delta.
  \end{align*}
  If we take \(j_m = j,\) \(n_m = n\) large enough also that \eqref{eq:c-accurate} 
  holds, then \(f_\infty\) satisfies \eqref{eq:f-indistinguishable}, and our 
  claim is proved.
\end{proof}

\subsection{Proof of Part B of Theorem \ref{rv}} \label{lpt}

\begin{proof}
 
  Suppose such \(C_n\) exists for $R=2r$. Set \(f_0 = 1,\) and
  \[f_1 = 1 + B2^{-j(r+1/2)} \sum_{k \in \mathcal Z_j} \beta_{jk} \psi_{jk},\]
  for \(B > 0,\) \(j > j_0,\) and \(\beta_{jk} = \pm 1\) to be determined.  
  Having chosen \(B,\) we will pick \(j\) large enough that \(f_1 \ge 
  \tfrac12.\) Since \(\int f_1 = \langle f_1, 1 \rangle = 1,\) \(f_1\) is then 
  a density.

  Set \(\delta = \tfrac14(1 - 2\alpha).\) As \(f_0 \in \Sigma(R, 1),\) for 
  \(n\) and $L$ large we have:
  \begin{enumerate}
    \item \(\Pr_{f_0}\{f_0 \not \in C_n\} \le \alpha + \delta;\) and
    \item \(\Pr_{f_0}\{\abs{C_n} \ge Ln^{-R/(2R+1)}\} \le \delta.\)
  \end{enumerate}
  Setting \(T_n = 1(\exists\, f \in C_n : \norm{f - f_0}_2 \ge 
  2Ln^{-R/(2R+1)}),\) we then have
  \[\Pr_{f_0}\{T_n = 1\} \le \alpha + 2\delta,\]
  as in the proof of Theorem \ref{imp}.

  For a constant \(C = C(\delta) > 0\) to be determined, set \(B = 
  (3L)^{2R+1}C^{-R}.\) Allowing \(j \to \infty,\) set \(n \sim 
  CB^{-2}2^{j(R+1/2)}.\) Then
  \[\norm{f_1 - f_0}_2 = B2^{-jr} \simeq 3Ln^{-R/(2R+1)},\]
  so for \(j\) large, \(\norm{f_1 - f_0}_2 \ge 2Ln^{-R/(2R+1)}.\) Arguing as 
  in the proof of Theorem \ref{imp}, the density \(Z\) of \(f_1\) w.r.t.\ \(f_0\) has second moment
  \begin{align*}
    E_{f_0}[Z^2] &\le \cosh(nB^22^{1-j(2r+1)})^{2^j}\\
    &= \cosh(C2^{1-j/2}(1 + o(1)))^{2^j}\\
    &= (1 + C^22^{2-j}(1 + o(1)))^{2^j}\\
    &\le \exp(4C^2(1 + o(1)))\\
    &\le 1 + \delta^2,
  \end{align*}
  for \(C(\delta)\) small, \(j\) large. Hence
  \[\Pr_{f_0}\{T_n=1\} + \max_\beta \Pr_{f_1}\{T_n = 0\} \ge 1 - \delta.\]
  and for all \(j\) (and \(n\)) large enough, we obtain, for suitable 
  \(\beta,\)
  \[\Pr_{f_1}\{f_1 \in C_n\} \le \Pr_{f_1}\{T_n = 1\} \le \alpha + 3\delta = 1 - 
  \alpha - \delta.\]
  Since \(f_1 \in \Sigma(r, B)\) for all $n, \beta_{jk}$ this contradicts the 
  definition of \(C_n.\)
\end{proof}

\bigskip

\textbf{Acknowledgement.} The authors are very grateful to two anonymous referees for a careful reading of a preliminary manuscript that led to several substantial improvements.

\bibliographystyle{plain}
\bibliography{band-bib}

\end{document}